\documentclass[twocolumn,letterpaper]{autart}
\usepackage{calc}
\let\theoremstyle\relax
\usepackage{amsmath,amssymb,amsthm}
\usepackage{tikz,pgfplots} 
\usepackage{natbib,enumerate}
\usepackage[short,nocomma]{optidef}
\usepackage{comment}
\usepackage{booktabs}

\usetikzlibrary{patterns}



\begingroup
    \makeatletter
    \@for\theoremstyle:=definition,remark,plain\do{%
        \expandafter\g@addto@macro\csname th@\theoremstyle\endcsname{%
            \addtolength\thm@preskip\parskip
            }%
        }
\endgroup

\theoremstyle{remark}

\theoremstyle{definition}

\newtheorem{assumption}{Assumption}

\theoremstyle{plain}

\newtheorem{theorem}{Theorem}

\newtheorem{definition}{Definition}
\newtheorem{lemma}{Lemma}
\newtheorem{remark}{Remark}

\allowdisplaybreaks[1]

\begin{document}

\begin{frontmatter}
\title{Feedback Nash equilibria for scalar N-player linear quadratic dynamic games\thanksref{footnoteinfo}}
\thanks[footnoteinfo]{The work of B. Nortmann has been partially supported by UKRI DTP 2019 grant no. EP/R513052/1.}
        
\author[Lon]{Benita Nortmann},
\ead{benita.nortmann15@imperial.ac.uk}
\author[Rome]{Mario Sassano}
\ead{mario.sassano@uniroma2.it}
\and
\author[Lon]{Thulasi Mylvaganam},
\ead{t.mylvaganam@imperial.ac.uk}

\address[Lon]{Department of Aeronautics, Imperial College London, London SW7 2AZ, UK }
\address[Rome]{Dipartimento di Ingegneria  Civile e Ingegneria Informatica, Universit{\`a} di Roma Tor Vergata, 00133 Roma, Italy}

\begin{abstract}  
Considering infinite-horizon, discrete-time, linear quadratic, $N$-player dynamic games with scalar dynamics, a graphical representation of feedback Nash equilibrium solutions is provided.
This representation is utilised to derive conditions for the number and properties of different feedback Nash equilibria a game may admit. The results are illustrated via a numerical example.
\end{abstract}

\end{frontmatter}

\section{Introduction} \label{sec:intro}
Dynamic game theory (\cite{BasarOlsder1998}) provides powerful tools to model dynamic interactions between strategic decision makers, with applications including, for instance, economics and ecology (\cite{jorgensen2007}), robotics (\cite{Li2019}), multi-agent systems (\cite{Mylvaganam2017,Cappello2021}), power systems (\cite{pratap2017}), and cyber-security (\cite{Basar2019}).
The \emph{Nash equilibrium} is a commonly considered solution concept and of natural interest in non-cooperative settings, since a unilateral deviation from the equilibrium by any player results in a worse outcome for the deviating player (\cite{BasarOlsder1998}). 
For infinite-horizon linear quadratic dynamic games, (linear) feedback Nash equilibrium (FNE) solutions, that is Nash equilibria involving linear static state-feedback strategies, are characterised via the solutions of coupled algebraic matrix equations. While reminiscent of the algebraic Riccati equations arising in linear quadratic optimal control, the coupled equations associated with FNE are generally difficult to solve and may admit multiple solutions with different outcomes (see, e.g. \cite{Starr1969,BasarOlsder1998,engwerda_2005}).
For the continuous-time setting, namely linear quadratic \emph{differential} games, the existence, number and properties of FNE solutions have been extensively studied (see, e.g. \cite{Papavassilopoulos1979,BasarOlsder1998,engwerda_2005,Possieri2015,ENGWERDA2016FNE}), however, the discrete-time equivalent has received significantly less attention (\cite{BasarOlsder1998,Pachter2010}). Additional product terms of the decision variables introduce further challenges specific to
the discrete-time case and make this case an interesting problem to study (see, \emph{e.g.} \cite{Monti2023, Nor-Mon:TAC24}).

With the aim of developing intuition regarding FNE solutions of discrete-time, linear quadratic, dynamic games, the focus of this note lies on a class of games in which the state and the players' inputs are scalar variables. For this class of games, a graphical representation of the coupled equations associated with FNE solutions is proposed. Utilising geometric arguments, conditions in terms of the system and cost parameters characterising the number and certain properties of FNE solutions a given game admits are derived.
The analysis can be considered a discrete-time counterpart to the study of scalar differential games presented in \cite{ENGWERDA2016FNE} and \cite[Chapter 8.4]{engwerda_2005}.
Despite similarities in the underlying constructions, a few noteworthy differences between the continuous-time and the discrete-time settings are observed and discussed in this paper.
Some preliminary results limited to games involving two players have appeared in \cite{scalarIFAC}. In this note, the results are generalised to the $N$-player case and a wider range of cost parameters are considered.

The remainder of the paper is organised as follows. Preliminaries on linear quadratic dynamic games are recalled in Section~\ref{sec:prelim}. A graphical interpretation of FNE solutions is proposed in Section~\ref{sec:graphical_int}. In Section~\ref{sec:FNE_cond}, this is utilised to derive conditions for the number and properties of FNE solutions a game may admit. The results are illustrated via a numerical example in Section~\ref{sec:example} and concluding remarks are provided in Section~\ref{sec:conclusion}.

\section{Preliminaries} \label{sec:prelim}
Consider a discrete-time linear time-invariant system influenced by the control actions $u_i \in \mathbb{R}$, $i = 1, \ldots, N$, of $N \in \mathbb{N}$ players, described by the scalar dynamics 
\begin{equation}
    x(k+1) = ax(k) + \sum_{i=1}^N b_i u_i(k),
    \label{eq:system}
\end{equation}
with $x(0) = x_0$, where $x \in \mathbb{R}$ denotes the system state and $a \in \mathbb{R}$, $b_i \in \mathbb{R} \setminus \{ 0 \} $, $i = 1, \ldots,N$, are constant system parameters\footnote{If $b_i = 0$, for any $i = 1, \ldots, N$, then $u_i$ does not influence the dynamics \eqref{eq:system} and player $i$ can be disregarded. Hence, this case is excluded without loss of generality.}. Let the quadratic cost functional 
\begin{equation}
    J_i(x_0, u_1, \ldots, u_N) = \sum_{k=0}^\infty \left( q_i x(k)^2 +  r_i u_i(k)^2 \right), 
    \label{eq:cost}
\end{equation}
with $q_i \in \mathbb{R}$, 
$r_i \in \mathbb{R}_{>0}$, be associated with player $i$, $i = 1, \ldots, N$. The dynamics \eqref{eq:system} and the cost functionals \eqref{eq:cost}, $i = 1, \ldots, N$, constitute a linear quadratic dynamic game.  
\begin{definition} \label{def:NE}
An admissible\footnote{A set of feedback strategies $\{u_1, \ldots, u_N\}$, with $u_i = k_i x$, $i = 1, ..., N$, is \emph{admissible} if it renders the zero equilibrium of system \eqref{eq:system} asymptotically stable.}
set of strategies $\{u_1^\star, \ldots, u_N^\star \}$ constitutes a Nash equilibrium solution of the game \eqref{eq:system}, \eqref{eq:cost}, $i = 1, \ldots, N$, if
\begin{equation}
    J_i^\star = J_i(x_0, u_i^\star, u_{-i}^\star) \leq J_i(x_0, u_i, u_{-i}^\star),
    \label{eq:NE}
\end{equation}
holds for all admissible $\{u_i, u_{-i}^\star\}$, for $i = 1, \ldots,N$, where 
\begin{equation*}
    u_{-i} = \{ u_1, \ldots, u_{i-1}, u_{i+1}, \ldots, u_N \}.
\end{equation*}
The strategy $u_i^\star$ is referred to as a Nash equilibrium strategy of player $i$, $i = 1, \ldots, N$, and the set $\{ J_1^\star, \ldots, J_N^\star\}$ is the Nash equilibrium outcome.
\end{definition}
Focusing on linear state-feedback strategies, Nash equilibrium solutions can be characterised via the stabilising solutions of a set of coupled algebraic equations. This characterisation of FNE solutions is recalled in the following statement.
\begin{theorem}
\label{th:coupled_algebraic_eq}
    Consider the game \eqref{eq:system}, \eqref{eq:cost}, $i = 1, \ldots, N$.
    The set of strategies $\{u_1^\star, \ldots, u_N^\star \}$, where
    \begin{equation}
        \label{eq:u_FNE}
        u_i^\star(k) = k_i x(k)\,,
    \end{equation}
    \noindent for $i = 1,\ldots, N$, constitutes a FNE of the game if and only if 
    \begin{equation}
        \left|a_{cl} \right| < 1,
        \label{eq:stability_cond}
    \end{equation}
    where $a_{cl} := a + \sum_{j=1}^N b_j k_j$, and there exist $p_i \in \mathbb{R}$, for $i = 1,\ldots, N$, satisfying the set of equations
    \begin{subequations}
        \label{eq:ACE}
        \begin{align}
            0 &= \left(a_{cl}^2 - 1 \right) p_i + q_i + k_i^2 r_i\,,
            \label{eq:ACE_1} \\[.5em]
            0 &= \left(r_i + b_i^2 p_i \right) k_i + b_i p_i \left( a + \sum_{j = 1, j \neq i}^N b_j k_j \right), 
            \label{eq:ACE_2}
        \end{align} 
        and such that
        \begin{equation}
            \left(r_i + b_i^2 p_i \right) > 0,
            \label{eq:2nd_order_condition}
        \end{equation}
    \end{subequations}
    for $i =1,\ldots, N$. The cost incurred by player $i$ starting from initial condition $x(0) = x_0$ is $J^\star_i =p_i x_0^2$.
\end{theorem}
\vspace{-0.45cm}
\begin{proof}
The proof follows from dynamic programming (\cite{Bellman1957}) via analogous arguments as used in the proof of \cite[Theorem 3.2]{Monti2023}, which considers two-player games with $q_i \geq 0$, for $i = 1,2$. In the more general case considered herein the condition \eqref{eq:2nd_order_condition} ensures that the minimisation problem for player $i$ is well-posed, for $i=1, \ldots, N$.
\end{proof}
Even in the scalar case considered herein, the coupled algebraic equations \eqref{eq:ACE}, $i = 1, \ldots, N$, are generally challenging to solve and there may be multiple stabilising solutions, as is the case for the continuous-time counterpart, see \emph{e.g.} \cite{engwerda_2005}.
In the following section, we use geometric arguments to derive conditions for the number and properties of FNE of the game \eqref{eq:system}, \eqref{eq:cost}, $i = 1,\ldots,N$.
\begin{remark}
\label{re:compare_c.t.}
It is interesting to note that contrary to their continuous-time counterpart (the coupled algebraic Riccati equations, see \emph{e.g.
\cite[Equation (3), $i = 1, \ldots,N$]{ENGWERDA2016FNE}}), \eqref{eq:ACE}, $i = 1, \ldots, N$, are \emph{not} quadratic equations. In fact, they are \emph{cubic} in the decision variables. Moreover, while in the continuous-time case the equilibrium gain of player $i$ depends on the gains of the other players only implicitly via the coupling of the equivalent to \eqref{eq:ACE_1}, namely \cite[Equation (3)]{ENGWERDA2016FNE},
it is evident from \eqref{eq:ACE_2} that in the discrete-time case $k_i$ explicitly depends on the equilibrium gains of all other players.
Finally, while in the continuous-time setting the condition $r_i > 0$ is sufficient to ensure the solutions of the equivalent of \eqref{eq:ACE_1}, \eqref{eq:ACE_2} correspond to minimising \eqref{eq:cost} for player $i$ (and hence any stabilising solution of \cite[Equation (3), $i = 1, \ldots,N$]{ENGWERDA2016FNE} constitutes a FNE), the additional condition \eqref{eq:2nd_order_condition} is needed in the discrete-time case, for $i = 1, \ldots, N$.
These differences introduce additional challenges in finding FNE solutions compared to the continuous-time setting, and make the discrete-time case interesting to study.
\end{remark}

\section{Graphical interpretation of FNE} \label{sec:graphical_int}
With the aim of deriving conditions for the existence of no, a unique or multiple FNE solutions of the game \eqref{eq:system},
\eqref{eq:cost}, $i = 1, \ldots, N$, we introduce a reformulation of Theorem~\ref{th:coupled_algebraic_eq}. To this end, consider the following assumptions. 
\begin{assumption} \label{as:player_order}
Let $\sigma_i := \frac{b_i^2q_i}{r_i}$, for $i = 1, \ldots, N$. The players are ordered such that $\sigma_1 \geq \sigma_2 \geq \ldots \geq \sigma_N$.
\end{assumption}
\begin{assumption}
    \label{as:enure_2nd_order_conditions}
    Let $\sigma_N > -1$. 
\end{assumption}
\begin{lemma} 
    \label{le:eliminate_2nd_order_conditions}
    Consider the game \eqref{eq:system}, \eqref{eq:cost}, $i = 1, \ldots, N$. If Assumptions~\ref{as:player_order} and \ref{as:enure_2nd_order_conditions} hold, then \eqref{eq:2nd_order_condition} holds for any solution of \eqref{eq:stability_cond}, \eqref{eq:ACE_1}, \eqref{eq:ACE_2}, for $i = 1, \ldots, N$.
\end{lemma}
\begin{proof}
    For \eqref{eq:2nd_order_condition} to hold we require $p_i > -\frac{r_i}{b_i^2}$. Consider first the case $a_{cl} \neq 0$. Combining \eqref{eq:ACE_1} and \eqref{eq:ACE_2}, and solving \eqref{eq:ACE_1} for $p_i$ gives that $p_i > -\frac{r_i}{b_i^2}$ 
    is equivalent to 
    \begin{equation}
        \frac{b_i k_i}{a_{cl}} < 1,
        \label{eq:alt_2nd_order_condition}
    \end{equation}
    for $i = 1, \ldots, N$.
    Combining again \eqref{eq:ACE_1} and \eqref{eq:ACE_2} by solving \eqref{eq:ACE_1} for $p_i$, substituting this into \eqref{eq:ACE_2} and multiplying by $\frac{b_i}{r_i}$ yields $0 = a_{cl}(b_i k_i)^2 + (1-a_{cl}^2) b_i k_i + a_{cl} \sigma_i$. Then, noting that $a_{cl} = a_i+ b_i k_i$, with $a_i := a + \sum_{j = 1, j \neq i}^N b_j k_j$, gives the condition
    $$ 0 = \left( b_i k_i\right)^2 + \left(a_i - \frac{(\sigma_i+1)}{a_i} \right) \left( b_i k_i\right) - \sigma_i,$$ for $i = 1, \ldots, N$.
    Solving this quadratic equation gives
    $b_i k_i = - a_i +\gamma_i \pm \sqrt{\gamma_i^2-1} $ and hence $a_{cl} = \gamma_i \pm \sqrt{\gamma_i^2-1}$, where $\gamma_i := \frac{1}{2} \left( a_i + \frac{\sigma_{i}+1}{a_i}\right)$. Thus, \eqref{eq:alt_2nd_order_condition} is in turn equivalent to
    \begin{equation*}
        1-\frac{a_i}{\gamma_i \pm \sqrt{\gamma_i^2-1}} < 1,
    \end{equation*}
    which holds if $a_i$ and $\gamma_i$ have the same sign.
    This in turn holds true if $\sigma_i > -1$.
    If $a_{cl}=0$, \eqref{eq:ACE_1} and \eqref{eq:ACE_2} imply $p_i = q_i$. Hence, \eqref{eq:2nd_order_condition} holds if $q_i > -\frac{r_i}{b_i^2}$, which is again equivalent to $\sigma_i > -1$.
    By Assumption~\ref{as:player_order}, Assumption~\ref{as:enure_2nd_order_conditions} implies $\sigma_i > -1$, for $i = 1, \ldots,N$. 
\end{proof}
\begin{remark}
    \label{re:assumptions}
    While Assumption~\ref{as:player_order} can be introduced without loss of generality, Assumption~\ref{as:enure_2nd_order_conditions} only depends on system and cost parameters and can hence be verified \emph{prior} to the computation of solutions. The
    relevance of Assumption~\ref{as:enure_2nd_order_conditions} is highlighted in Lemma~\ref{le:eliminate_2nd_order_conditions}. Note that the results presented in the remainder of the paper are still relevant if Assumption~\ref{as:enure_2nd_order_conditions} does not hold. However, in this case they concern any solutions of \eqref{eq:stability_cond}, \eqref{eq:ACE_1}, \eqref{eq:ACE_2}, $i=1, \ldots,N$. Hence, the condition \eqref{eq:2nd_order_condition}, or alternatively \eqref{eq:alt_2nd_order_condition}, needs to be checked to ensure such a solution corresponds to a FNE.
\end{remark}

\begin{lemma} \label{le:aux_fun_desc}
Consider the game \eqref{eq:system}, \eqref{eq:cost}, $i = 1, \ldots, N$. Let Assumptions~\ref{as:player_order} and \ref{as:enure_2nd_order_conditions} hold and consider the function
 \begin{align}
    \label{eq:hat_f}
    \hat{f}(\xi) = 
    \begin{cases}
        -\xi-\sqrt{\xi^2+1} & \text{if } \xi < 0, \\
        -\xi+\sqrt{\xi^2+1} & \text{if } \xi > 0.
    \end{cases}
\end{align}
The set of strategies $\{u_1^\star, \ldots, u_N^\star \}$, where $u_i^\star(k)$ is given by \eqref{eq:u_FNE} with 
\begin{equation}
    k_i = \frac{-\xi - t_i \sqrt{\xi^2- \sigma_i}}{b_i}, 
    \label{eq:k_lem_aux}
\end{equation}
for $i = 1,\ldots,N$, constitutes a FNE of the game with $a_{cl} \neq 0$, if and only if there exist $t_i \in \left\{-1,1 \right\}$, for $i= 1,\ldots,N$, and $\xi \in \mathbb{R} \setminus \{ 0 \}$, satisfying
\begin{align}
    \label{eq:aux_sol}
    a = \hat{f}(\xi) + N \xi + t_1\sqrt{\xi^2 - \sigma_1} + \ldots + t_N\sqrt{\xi^2-\sigma_N}.
\end{align}
\end{lemma}
\begin{proof}
    By Theorem~\ref{th:coupled_algebraic_eq}, FNE solutions of the game are characterised by the stabilising solutions of \eqref{eq:ACE}, $i = 1,\ldots,N$. The proof lies in showing that solving \eqref{eq:aux_sol} is equivalent to solving \eqref{eq:stability_cond}, \eqref{eq:ACE}, $i = 1,\ldots,N$. Eliminating $p_i$ in \eqref{eq:ACE} (namely solving \eqref{eq:ACE_1} for $p_i$ and substituting this into \eqref{eq:ACE_2}) gives the condition
    \begin{equation}
        \label{eq:reform_ARE}
        0 = \frac{b_i}{2}k_i^2 + \xi k_i + \frac{\sigma_i}{2b_i},
    \end{equation}
    with $\xi := \frac{1}{2} \left( \frac{1}{a_{cl}} - a_{cl}\right)$. The equation \eqref{eq:reform_ARE} admits the solutions \eqref{eq:k_lem_aux}, $t_i \in \left\{-1,1 \right\}$, for $i = 1, \ldots,N$. Hence, $a_{cl}$ can be written as
    \begin{align}
        \label{eq:aux_acl}
        \begin{split}
        a_{cl} = -\xi \pm \sqrt{\xi^2 + 1} &= a + \sum_{j=1}^N b_j k_j, \\
        &= a - N\xi - t_1 \sqrt{\xi^2 - \sigma_1} - \ldots \\
        & \hspace{2.3cm} - t_N \sqrt{\xi^2 - \sigma_N}.
        \end{split}
    \end{align}
    By Definition~\ref{def:NE}, we are interested in solutions \eqref{eq:k_lem_aux}, $t_i \in \left\{-1,1 \right\}$, such that \eqref{eq:stability_cond} holds. Hence, there is a one-to-one correspondence between $\xi$ and $a_{cl}$ given by $a_{cl} = \hat f(\xi)$ as defined in \eqref{eq:hat_f}. Substituting this into \eqref{eq:aux_acl} gives \eqref{eq:aux_sol}. Thus, any solution to \eqref{eq:stability_cond}, \eqref{eq:ACE}, $i = 1, \ldots,N$, is such that \eqref{eq:aux_sol} holds. 
    Conversely, let $\{t_1,\ldots,t_N \}$ and $\xi$ be a solution to \eqref{eq:aux_sol}. Note that this solution is such that \eqref{eq:reform_ARE} holds with $k_i$ as in \eqref{eq:k_lem_aux}, for $i = 1,\ldots,N$. Via $a_{cl} = \hat f(\xi)$ as defined in \eqref{eq:hat_f}, \eqref{eq:stability_cond} holds and \eqref{eq:reform_ARE} implies that \eqref{eq:ACE_2} holds with $p_i = \frac{q_i + k_i^2 r_i}{1-a_{cl}^2}$, which is the unique solution of the Lyapunov equation \eqref{eq:ACE_1} for fixed $a_{cl}$, for $i=1,\ldots,N$. By Lemma \ref{le:eliminate_2nd_order_conditions}, Assumption~\ref{as:enure_2nd_order_conditions} ensures that \eqref{eq:2nd_order_condition} holds, for $i=1,\ldots,N$. 
    Hence, \eqref{eq:aux_sol} implies \eqref{eq:stability_cond}, \eqref{eq:ACE}, $i = 1,\ldots,N$. 
\end{proof}
\begin{remark}
    The result
    of Lemma~\ref{le:aux_fun_desc} 
    introduces the assumption $a_{cl} \neq 0$. Note that if $a_{cl} = 0$, then Theorem~\ref{th:coupled_algebraic_eq}, in particular \eqref{eq:ACE_2}, implies $k_i = 0$, for $i=1,\ldots,N$, and hence $a=0$. Thus, the assumption $a_{cl} \neq 0$ is only restrictive in the special case in which $a = 0$. In this case, a set of FNE strategies is given by \eqref{eq:u_FNE}, with $k_i = 0$, for $i=1,\ldots,N$. 
    However, this trivial solution cannot be found via the result of Lemma~\ref{le:aux_fun_desc}.
    All FNE solutions of the game \eqref{eq:system}, \eqref{eq:cost}, $i = 1, \ldots, N$, with $a=0$ are hence given by the solutions satisfying the conditions of Lemma~\ref{le:aux_fun_desc} (if there are any) and, in addition, the solution \eqref{eq:u_FNE}, with $k_i = 0$, 
    $i=1,\ldots,N$. 
    \label{re:a=0}
\end{remark}
In Lemma~\ref{le:aux_fun_desc} FNE solutions of the game  \eqref{eq:system}, \eqref{eq:cost}, $i = 1, \ldots, N$, are characterised via the condition \eqref{eq:aux_sol}.
Consider the auxiliary functions
\begin{equation}
    f_\ell(\xi) = \hat f(\xi) + N \xi + \tau_{\ell, 1} \sqrt{\xi^2 - \sigma_1} + \ldots + \tau_{\ell, N}\sqrt{\xi^2-\sigma_N},
    \label{eq:aux_fun}%
    \smallskip
\end{equation}
for $\ell = 1, \ldots, L$, where $L = 2^N$ and $\tau_\ell = \left( \tau_{\ell, 1}, \ldots, \tau_{\ell, N} \right)$ is an $N$-tuple over the set $\{ -1, 1\}$. 
The functions $f_\ell(\xi)$, $\ell = 1, \ldots, L$, in \eqref{eq:aux_fun} capture all possible combinations of the values which $t_i$, for $i = 1, \ldots, N$, can take in \eqref{eq:aux_sol}.
Hence, by Lemma~\ref{le:aux_fun_desc}, FNE solutions of the game \eqref{eq:system}, \eqref{eq:cost}, $i = 1, \ldots, N$, are represented graphically by the intersections of the functions $f_\ell(\xi)$, $\ell = 1, \ldots, L$, as defined in \eqref{eq:aux_fun}, with the horizontal line at level $a$.
\begin{remark}
    \label{re:multiple_intersections}
    If the horizontal line at level $a$ intersects multiple auxiliary functions $f_\ell(\xi)$, $\ell = 1, \ldots, L$, as defined in \eqref{eq:aux_fun}, in a point in which they coincide, then this intersection point generally corresponds to multiple \emph{distinct FNE solutions} of the game \eqref{eq:system}, \eqref{eq:cost}, $i = 1, \ldots, N$, resulting in the \emph{same closed-loop dynamics} $a_{cl}$. 
    However, consider the special case of two functions $f_\ell(\xi)$ and $f_w(\xi)$, for $\ell = 1, \ldots, L$, $w = 1, \ldots, L$, $\ell \neq w$, which coincide in the point $(\bar \xi, \bar f)$ with $\bar \xi = \pm \sqrt{\sigma_j}$, $j = 1, \ldots, N$, and $\bar f = f_\ell(\bar \xi) = f_w(\bar \xi)$. Moreover, $f_\ell(\xi)$ and $f_w(\xi)$ are such that $\tau_{\ell,i} = \tau_{w,i}$, for all $i = 1, \ldots,N$, $i \neq j$, or for any $i \neq l$, $i \neq j$, if the game is such that any $\sigma_l = \sigma_j$, $l = 1, \ldots, N$, $l\neq j$.
    If the line at level $a$ intersects these two functions in $(\bar \xi, \bar f)$, then there exists only one corresponding set of gains $\{k_1, \ldots, k_N\}$. This set of gains is given by $k_j = -\frac{\bar \xi}{b_j}$, $k_l = -\frac{\bar \xi}{b_l}$ and $k_i$ as defined in \eqref{eq:k_lem_aux} with $t_i = \tau_{\ell,i} = \tau_{w,i}$, for $i \neq j$, $i \neq l$.
    Hence, in this special case the intersection point $(\bar \xi, \bar f)$ corresponds to one FNE solution rather than two distinct FNE solutions.
\end{remark}
With the aim of utilising this graphical representation to characterise the possible number and properties of FNE, we first take a closer look at the functions \eqref{eq:aux_fun}, $\ell = 1, \ldots, L$.
Let $\mathcal{T} = \left\{\tau_1, \ldots, \tau_{L} \right\}$ denote the set of all $N$-tuples over $\{ -1, 1\}$ and consider the function $T_\ell : \{1, \ldots, N \} \rightarrow \{ \tau_{\ell, 1}, \ldots, \tau_{\ell, N} \}$, defined by $T_\ell(i) = \tau_{\ell, i}$, for $i = 1, \ldots, N$.
Consider in particular the functions $f_l$, $l = 1,2,3, L-2, L-1, L$, and define the N-tuples $\tau_l \in \mathcal{T}$, for $l = 1,2,3,L-2,L-1,L$, such that 
\begin{align*}
     &T_1(i) = -1, &&i = 1, \ldots, N,\\
     &T_2(1) = 1, T_2(i) = -1, &&i = 2, \ldots, N ,\\
     &T_3(2) = 1, T_3(i) = -1, &&i = 1, \ldots,N, i \neq 2, \\
     &T_{L-2}(2) = -1,  T_{L-2}(i) = 1, &&i = 1, \ldots,N, i \neq 2, \\
     &T_{L-1}(1) = -1, T_{L-1}(i) = 1, &&i = 2, \ldots, N ,\\
     &T_L(i) = 1, &&i = 1, \ldots, N.
\end{align*}
With this notation in place, let us highlight some properties of the auxiliary functions \eqref{eq:aux_fun}, $\ell = 1, \ldots, L$, in the following statement\footnote{Note that if $N \geq 3$ and hence $L > 6$, the $N$-tuples $\tau_{\ell}$ for $\ell = 4, \ldots, L-3$, and hence the functions $f_\ell$, $\ell = 4, \ldots, L-3$, are not defined explicitly, since they are not relevant for the following analysis.}.
\begin{lemma}
    \label{le:function_properties}
    Consider the functions \eqref{eq:aux_fun},  $\ell = 1, \ldots, L$, and let Assumption~\ref{as:player_order} hold. 
    Then,
    \begin{enumerate}[i.]
        \item $f_1 \leq f_2 \leq f_3 \leq f_\ell \leq f_{L-2} \leq f_{L-1} \leq f_{L}$, for any $\ell = 4, \ldots, L-3$.
        \label{item:aux_fun_ordering}
        \item $\displaystyle \lim_{\xi \rightarrow - \infty} \left| f_\ell (\xi)  - \left( N - \sum_{j = 1}^N \tau_{\ell,j} \right) \xi \right| = 0$, and \\ 
        $\displaystyle \text{ }\lim_{\xi \rightarrow \infty} \hspace{0.05cm} \left| f_\ell (\xi) - \left( N + \sum_{j = 1}^N \tau_{\ell,j} \right) \xi \right| = 0$.
        \label{item:aux_fun_asymptotes}
        \item $f_\ell(\xi)$ is defined over the real numbers for
        \begin{enumerate}
            \item $\xi \neq 0$, for $\ell = 1, \ldots, L$, if $\sigma_1 \leq 0$. 
            \item $\xi \leq -\sqrt{\sigma_1}$ and $\xi \geq \sqrt{\sigma_1}$, for $\ell = 1, \ldots, L$, if $\sigma_1 > 0$ and $\sigma_1 \neq \sigma_2$.
            \item $\xi \leq -\sqrt{\sigma_1}$ and $\xi \geq \sqrt{\sigma_1}$, for $\ell = 1,L$, and \\ for $ \begin{cases}
            \xi \neq 0 & \text{if } \sigma_3 \leq 0,\\
            \xi \leq -\sqrt{\sigma_3} \text{ and } \xi \geq \sqrt{\sigma_3}& \text{if } \sigma_3 > 0,   
            \end{cases}$ \\
            for $\ell = 2,3,L-2, L-1$, if $\sigma_1 > 0$ and $\sigma_1 = \sigma_2$.
        \end{enumerate}
        \label{item:aux_fun_definition}
        \item if $\sigma_1 > 0$ and $\sigma_1 = \sigma_2$, and there exist $\sigma_j > 0$ such that $\sigma_j \neq \sigma_i$, for $i = 1, \ldots,N$, $j= 3, \ldots,N$, $i \neq j$, then let $\bar \sigma = \max_j (\sigma_j)$. No function $f_{\ell}(\xi)$, $\ell =  1, \ldots, L$, is defined over the real numbers for $-\sqrt{\bar \sigma} < \xi < \sqrt{\bar \sigma}$.
        \label{item:aux_fun_not_defined}
        \item if $f_\ell(\xi)$ is defined for $\xi \neq 0$, then \\
        $\displaystyle \lim_{\xi \rightarrow 0^{-}} f_\ell(\xi) = -1 + \sum_{i=1}^N \tau_{\ell,i} \sqrt{-\sigma_i} = \bar a_{\ell}^{-},$ and \\
        $\displaystyle \lim_{\xi \rightarrow 0^{+}} f_\ell(\xi) = 1 + \sum_{i=1}^N \tau_{\ell,i} \sqrt{-\sigma_i} = \bar a_{\ell}^{+}.$
        \label{item:aux_fun_xi=0}
        \item if for all $\xi \neq 0$ \\
        $\displaystyle (N-1) + \frac{|\xi|}{\sqrt{\xi^2+1}} \neq \sum_{i=1}^N \frac{|\xi|}{\sqrt{\xi^2-\sigma_i}},$\\
        then $f_L(\xi)$ is strictly monotone for $\xi < 0$ and $f_1(\xi)$ is strictly monotone for $\xi > 0$.
        \label{item:aux_fun_monotone}
    \end{enumerate}
\end{lemma}
\begin{proof}
    Recall that by Assumption~\ref{as:player_order} $\sigma_1 \geq \sigma_2 \geq \ldots \geq \sigma_N$. 
    Item~\emph{\ref{item:aux_fun_ordering}.} follows from the definition of $T_\ell$, for $l = 1,2,3,L-2,L-1,L$, and Assumption~\ref{as:player_order}. \\
    Item~\emph{\ref{item:aux_fun_asymptotes}.} follows from \eqref{eq:aux_fun} by noting that \sloppy
    $\lim_{\xi \rightarrow \infty} |\sqrt{\xi^2+c} - \xi | = 0 , \forall c \in \mathbb{R}, |c| < \infty.$ \\
    Item~\emph{\ref{item:aux_fun_definition}.} follows by noting that \eqref{eq:aux_fun} is defined over the real numbers if $\xi^2 - \sigma_i \geq 0$, for $i = 1, \ldots, N$, and utilising Assumption~\ref{as:player_order} and the definition of $T_\ell$, $l = 1,2,3,L-2,L-1,L$. \\
    Similarly, to demonstrate item~\emph{\ref{item:aux_fun_not_defined}.} note that if $\sigma_1 = \sigma_2$, then the corresponding terms in \eqref{eq:aux_fun}, $\ell = 1, \ldots, L$, may cancel out, which affects the interval over which some functions are defined. 
    If any $\sigma_j > 0$, $j = 3, \ldots, N$, such that $\sigma_j \neq \sigma_i$, for $i = 1, \ldots,N$, $i \neq j$, and $\bar \sigma = \max_j (\sigma_j)$, then the corresponding terms do not cancel in any functions. Utilising Assumption~\ref{as:player_order}, the functions in which all terms corresponding to repeated $\sigma_i$, $i = 1, \ldots,N$, cancel out are defined over the real numbers if $\xi^2 - \bar \sigma \geq 0$. \\
    Item~\emph{\ref{item:aux_fun_xi=0}.} follows from \eqref{eq:aux_fun} by noting that \eqref{eq:hat_f} is such that
    $\lim_{\xi \rightarrow 0^{-}} \hat f (\xi) = -1 $, and $\lim_{\xi \rightarrow 0^{+}} \hat f (\xi) = 1 $. \\
    Finally, item~\emph{\ref{item:aux_fun_monotone}.} is shown by noting that $\frac{d }{d \xi} f_\ell (\xi) = (N-1) \pm \frac{\xi}{\sqrt{\xi^2+1}} + \sum_{i=1}^N \tau_{\ell,i} \frac{\xi}{\sqrt{\xi^2-\sigma_i}}$,
    and using the definition of $T_\ell$, $\ell = 1,L$.
\end{proof}

\section{Conditions for existence, number and properties of FNE} \label{sec:FNE_cond}
Via the result of Lemma~\ref{le:aux_fun_desc}, FNE solutions of the game \eqref{eq:system}, \eqref{eq:cost}, $i = 1, \ldots,N$, can be represented graphically by the intersections of the functions \eqref{eq:aux_fun}, $\ell = 1, \ldots, L$, with the horizontal line at level $a$. 
The properties of the functions \eqref{eq:aux_fun}, $\ell = 1, \ldots, L$, highlighted in Lemma~\ref{le:function_properties} allow us to characterise the possible number and location of such intersection points, and hence to derive conditions for the existence of a FNE, as well as the number and properties of different FNE solutions of the game, in the following statement.
\begin{theorem}
\label{th:number&properties}
Consider the game \eqref{eq:system}, \eqref{eq:cost}, $i = 1, \ldots,N$, and let Assumptions~\ref{as:player_order} and \ref{as:enure_2nd_order_conditions} hold. Then,
\begin{enumerate}[i.]
    \item if system \eqref{eq:system} is open-loop unstable with a fast rate of divergence, i.e. $|a| \gg 1$, then there exist $2^N-1$ FNE solutions.
    \label{item:nr_FNE_a>>1}
    \item if $\sigma_1 > 0$, then there always exists at least one FNE solution for any value of $a$. 
    \label{item:at_least_1_FNE}
    \item if $\sigma_1 > 0$ and there exist $\sigma_j > 0$ such that $\sigma_j \neq \sigma_i$, for $i = 1, \ldots,N$, $j= 1, \ldots,N$, $j \neq i$, and $\bar \sigma = \max_j (\sigma_j)$, then the system \eqref{eq:system} in closed loop with any FNE solution is such that $\left| a_{cl}^\star \right| \leq \left|\sqrt{\bar \sigma} - \sqrt{\bar \sigma+1} \right| $.
    \label{item:acl_bound}
    \item if $\sigma_N \geq 0$ and the system \eqref{eq:system} is open-loop stable, i.e. $|a| < 1$, then there exists a unique FNE solution.
    \label{item:unique_FNE_stable_a}
    \item if $\sigma_N < 0$, but $\sigma_1 > 0$, $\sigma_1 > \sigma_2$ and the system \eqref{eq:system} is open-loop stable with a fast rate of convergence, i.e. $|a| \ll 1$, then there exists a unique FNE solution if 
    $\sum_{i=2}^N \sqrt{\sigma_1 - \sigma_i} < (N-1) \sqrt{\sigma_1} + \sqrt{\sigma_1 + 1}$ and $ (N-1) + \frac{| \xi |}{\sqrt{\xi^2 + 1}} < \sum_{i=1}^N \frac{|\xi|}{\sqrt{\xi^2-\sigma_i}}.$
    \label{item:unique_FNE_mixed}
    \item if $\sigma_1 \leq 0$, and the system \eqref{eq:system} is open-loop stable with a fast rate of convergence, i.e. $|a| \ll 1$, then there exists a unique FNE solution if \\
    $\left(\sum_{i=2}^N \sqrt{- \sigma_i} \right) - \sqrt{-\sigma_1} < 1$ and $ (N-1) + \frac{|\xi|}{\sqrt{\xi^2+1}} \neq \sum_{i=1}^N \frac{|\xi|}{\sqrt{\xi^2-\sigma_i}}.$
    \label{item:unique_FNE_negative}
    \item if $\sigma_1 = \ldots = \sigma_N = 0$ and $|a| = 1$ there exists no FNE solution.
    \label{item:no_FNE_symmetric}
\end{enumerate}
\end{theorem}
\begin{proof}
    By Lemma~\ref{le:aux_fun_desc} the intersection points of the functions~\eqref{eq:aux_fun}, $\ell = 1, \ldots, L$, with the horizontal line at level $a$ correspond to distinct FNE solutions of the game \eqref{eq:system}, \eqref{eq:cost}, $i = 1, \ldots,N$. The claims are hence shown by characterising the possible intersection points. \\
    Item~\emph{\ref{item:nr_FNE_a>>1}.} is a result of item~\emph{\ref{item:aux_fun_asymptotes}.} of Lemma~\ref{le:function_properties}. For large values of $a$, the line at level $a$ intersects once with $L-1 = 2^N-1$ of the $L$  functions \eqref{eq:aux_fun}, $\ell = 1, \ldots, L$. \\
    Item~\emph{\ref{item:at_least_1_FNE}.} follows from item~\emph{\ref{item:aux_fun_asymptotes}.} and item~\emph{\ref{item:aux_fun_definition}.(b),(c)} of Lemma~\ref{le:function_properties}. Note that by definition \eqref{eq:aux_fun}, $\ell = 1, \ldots, L$, are continuous for $\xi \neq 0$ and that $f_{L-1}(-\sqrt{\sigma_1}) = f_{L}(-\sqrt{\sigma_1})$ and $f_{1}(\sqrt{\sigma_1}) = f_{2}(\sqrt{\sigma_1})$. Hence, for any value of $a \neq 0$, the horizontal line at level $a$ intersects at least once with either of the functions $f_\ell(\xi)$, $\ell = 1,2,L-1, L$. The case $a = 0$ is discussed in Remark~\ref{re:a=0}. \\
    Item~\emph{\ref{item:acl_bound}.} follows from item~\emph{\ref{item:aux_fun_not_defined}.} of Lemma~\ref{le:function_properties}. If the stated conditions hold, any intersection points of  \eqref{eq:aux_fun}, $\ell = 1, \ldots,L$, with the horizontal line at level $a$ are such that $\xi^\star \leq -\sqrt{\bar \sigma}$ or $\xi^\star \geq \sqrt{\bar \sigma}$. The claim is shown by recalling that $a_{cl} = \hat f (\xi)$, with $\hat f (\xi)$ defined in \eqref{eq:hat_f}. \\
    Item~\emph{\ref{item:unique_FNE_stable_a}.} follows from  items~\emph{\ref{item:aux_fun_ordering}.}, \emph{\ref{item:aux_fun_asymptotes}.}, \emph{\ref{item:aux_fun_definition}.} and \emph{\ref{item:aux_fun_monotone}.} of Lemma~\ref{le:function_properties}. If $\sigma_N \geq 0$, then by Assumption~\ref{as:player_order} $\sigma_i \geq 0$, for $i = 1, \ldots, N$, and the line at level $a$ intersects only once with one of the functions \eqref{eq:aux_fun}, $\ell = 1, L$, for $\max_{\xi < 0}(f_{L-1}(\xi)) \leq -1 < a < 1 \leq \min_{\xi >0}(f_2(\xi))$.  \\
    Item~\emph{\ref{item:unique_FNE_mixed}.} follows from items~\emph{\ref{item:aux_fun_ordering}.}, \emph{\ref{item:aux_fun_asymptotes}.}, \emph{\ref{item:aux_fun_definition}.(b)} and \emph{\ref{item:aux_fun_monotone}.} of Lemma~\ref{le:function_properties}. If the stated conditions hold, then $f_L(-\sqrt{\sigma_1}) < f_1(\sqrt{\sigma_1})$ and $f_L(\xi)$ for $\xi < 0$ and $f_1(\xi)$ for $\xi > 0$ are strictly decreasing. Hence, the line at level $a$ intersects only once with one of the functions \eqref{eq:aux_fun}, $\ell = 1, L$, for very small values of $a$. \\
    Item~\emph{\ref{item:unique_FNE_negative}.} follows from items~\emph{\ref{item:aux_fun_ordering}.}, \emph{\ref{item:aux_fun_asymptotes}.}, \emph{\ref{item:aux_fun_definition}.(a)}, \emph{\ref{item:aux_fun_not_defined}.} and \emph{\ref{item:aux_fun_monotone}.} of Lemma~\ref{le:function_properties}. If the stated conditions hold, then $\bar a_{L-1}^{-} < \bar a_2^{+}$ and $f_L(\xi)$ for $\xi < 0$ and $f_1(\xi)$ for $\xi > 0$ are both either strictly decreasing or strictly increasing. Hence, the line at level $a$ intersects only once with one of the functions \eqref{eq:aux_fun}, $\ell = 1, L$, for very small values of $a$. \\
    Finally, item~\emph{\ref{item:no_FNE_symmetric}.} follows from  items~\emph{\ref{item:aux_fun_ordering}.}, \emph{\ref{item:aux_fun_asymptotes}.}, \emph{\ref{item:aux_fun_definition}.(a)} and \emph{\ref{item:aux_fun_xi=0}.} of Lemma~\ref{le:function_properties}. 
    Note that if $\sigma_i = 0$, for $i = 1, \ldots,N$, then $f_L(\xi) > -1$, $f_{L-1}(\xi) < -1$ for $\xi < 0$, $f_1(\xi) < 1$, $f_2(\xi)  > 1$ for $\xi > 0$, and $\bar a_\ell^{-} = -1$, $\bar a_\ell^{+} = 1$ for $\ell = 1, \ldots, L$.
    Hence, since \eqref{eq:aux_fun}, $\ell = 1, \ldots, L$, are not defined at $\xi = 0$, the horizontal lines at $a= 1$ and at $ a = -1$ do not intersect with any of the functions.
\end{proof}
\begin{remark}
    \label{re:result_compare_c.t.}
    It is interesting to compare the results of Theorem~\ref{th:number&properties} with the continuous-time counterpart. While analogous results to items~\textit{\ref{item:nr_FNE_a>>1}.} - \textit{\ref{item:unique_FNE_stable_a}.} and \textit{\ref{item:no_FNE_symmetric}.} can be derived for scalar linear quadratic \emph{differential games}, see \cite{ENGWERDA2016FNE}, note that for this class of continuous-time dynamic games there always exists a unique feedback Nash equilibrium for games involving open-loop stable systems with fast rate of convergence, irrespective of the the signs of $\sigma_i$, $i = 1, \ldots, N$, \cite[Theorem 3.1]{ENGWERDA2016FNE}. In the discrete-time setting, however, this only holds true for special cases if $\sigma_N < 0$, such as under the conditions of items~\textit{\ref{item:unique_FNE_mixed}.} and \textit{\ref{item:unique_FNE_negative}.} of Theorem~\ref{th:number&properties}.
\end{remark}


\section{Example} \label{sec:example}
\begin{figure}[t]
    \centering
    \includegraphics[width=\linewidth]{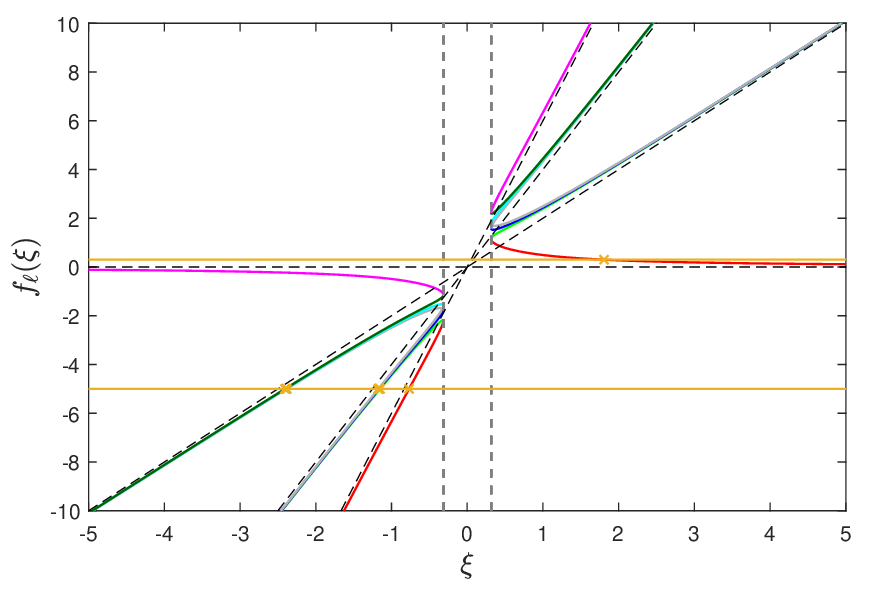}
    \caption{
    Plot of the auxiliary functions for $\tilde q_2 = 0.05$ and $\tilde q_3= 0$, $f_1(\xi)$ (red), $f_2(\xi)$ (green), $f_3(\xi)$ (blue), $f_\ell(\xi)$, $\ell = 4, \ 5$ (grey), $f_{6}(\xi)$ (cyan), $f_{7}(\xi)$ (dark green) and $f_{8}(\xi)$ (magenta), and their intersection points (yellow crosses) with the horizontal lines at $a = \tilde a$ (yellow) for $\tilde a = 0.3$ and $\tilde a = -5$. The black dashed lines indicate the linear asymptotes of $f_\ell(\xi)$, $\ell = 1, \ldots,8$, and the grey dotted lines indicate $\xi = \pm \sqrt{\sigma_1} = \pm 0.3162$.
    }
    \label{fig:q0}
\end{figure}

\begin{figure}[t]
    \centering
    \includegraphics[width=\linewidth]{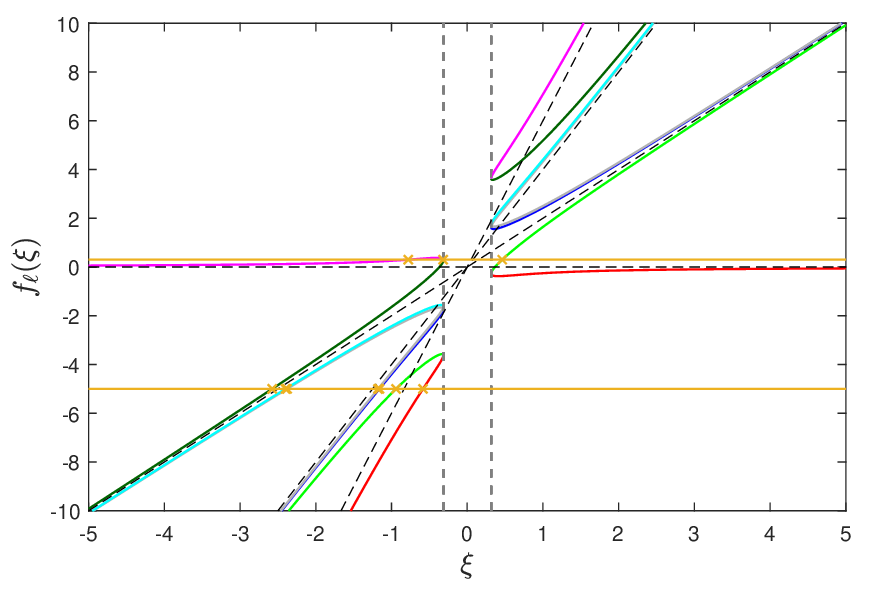}
    \caption{
    Plot of the auxiliary functions for $\tilde q_2 = -0.8$ and $\tilde q_3= -0.9$, $f_1(\xi)$ (red), $f_2(\xi)$ (green), $f_3(\xi)$ (blue), $f_\ell(\xi)$, $\ell = 4, \ 5$ (grey), $f_{6}(\xi)$ (cyan), $f_{7}(\xi)$ (dark green) and $f_{8}(\xi)$ (magenta), and their intersection points (yellow crosses) with the horizontal lines at $a = \tilde a$ (yellow) for $\tilde a = 0.3$ and $\tilde a = -5$. The black dashed lines indicate the linear asymptotes of $f_\ell(\xi)$, $\ell = 1, \ldots,8$, and the grey dotted lines indicate $\xi = \pm \sqrt{\sigma_1} = \pm 0.3162$.
    }
    \label{fig:qm5}
\end{figure}

To illustrate the presented results, consider the game defined by the dynamics \eqref{eq:system} with $N = 3$, $a = \tilde{a}$, and $b_1 = b_2 = b_3 = 1$, and the cost functionals \eqref{eq:cost}, $i = 1,2,3$, with $q_1 = 0.1$, $q_2 = \tilde q_2$, $q_3 = \tilde q_3$ and $r_1 = r_2 = r_3 = 1$. Note that $\sigma_i = q_i$, for $i = 1, 2, 3$. Firstly, let $\tilde q_2 = 0.05$ and $\tilde q_3= 0$. The corresponding auxiliary functions \eqref{eq:aux_fun}, $\ell = 1, \ldots, 8$, are plotted in Figure~\ref{fig:q0}. The horizontal yellow lines indicate the levels $\tilde a = 0.3$ and $\tilde a = -5$ and the yellow crosses highlight the intersection points between the horizontal lines and the functions \eqref{eq:aux_fun}, $\ell = 1, \ldots, 8$. By Lemma~\ref{le:aux_fun_desc} these intersection points correspond to FNE solutions of the game for the respective value of $\tilde a$.
For $\tilde a = 0.3$ there is only a single intersection point. In line with item~\emph{\ref{item:unique_FNE_stable_a}.} of Theorem~\ref{th:number&properties}, this single intersection point corresponds to a unique FNE solution of the game. 
For $\tilde a = -5$ there are $2^3-1 = 7$ intersection points, which in line with item~\emph{\ref{item:nr_FNE_a>>1}.} of Theorem~\ref{th:number&properties} correspond to $7$ FNE solutions.
Secondly, let $\tilde q_2 = -0.8$ and $\tilde q_3= -0.9$. The corresponding auxiliary functions \eqref{eq:aux_fun}, $\ell = 1, \ldots, 8$, are plotted in Figure~\ref{fig:qm5}. As above, the horizontal yellow lines indicate the levels $\tilde a = 0.3$ and $\tilde a = -5$ and the yellow crosses highlight the intersection points between the lines and the functions \eqref{eq:aux_fun}, $\ell = 1, \ldots, 8$. For $\tilde a = -5$ there are again $2^3-1 = 7$ intersection points. However, for $\tilde a = 0.3$ there are $3$ intersection points. By inspecting the functions $f_1(\ell)$, $f_2(\ell)$, $f_7(\ell)$ and $f_8(\ell)$ it is evident that there are more than $1$ intersection points for any $|a| \ll 1$. Note that the conditions $\sum_{i=2}^3 \sqrt{\sigma_1 -\sigma_i} < 2 \sqrt{\sigma_1} + \sqrt{\sigma_1+1}$ and $ 2 + \frac{| \xi |}{\sqrt{\xi^2 + 1}} < \sum_{i=1}^3 \frac{|\xi|}{\sqrt{\xi^2-\sigma_i}}$ are not satisfied for this example. Hence this observation is in line with item~\emph{\ref{item:unique_FNE_mixed}.} of Theorem~\ref{th:number&properties} and highlights the difference between discrete-time and continuous-time scalar linear quadratic dynamic games called to attention in Remark~\ref{re:result_compare_c.t.}. 
Note that for all described cases it holds that $|a_{cl}^\star| \leq |\sqrt{\sigma_1} - \sqrt{\sigma_1+1}| = 0.7326$ in line with item~\emph{\ref{item:acl_bound}.} of Theorem~\ref{th:number&properties}.

\vspace{-0.2cm}
\section{Conclusion} \label{sec:conclusion}
Considering infinite-horizon, discrete-time, linear quadratic, $N$-player dynamic games involving dynamics in which the state and the input of each player are scalar variables, a graphical representation of the conditions characterising FNE solutions is proposed. Via geometric arguments, this representation allows to derive conditions in terms of the system and cost parameters characterising the number and properties of solutions. The results are illustrated via a numerical example.
\vspace{-0.2cm}
\bibliographystyle{agsm}
\bibliography{biblio}

\end{document}